\documentclass[a4,16pt]{article} 
\usepackage{fontenc}
\usepackage{amsmath}
\usepackage{amssymb}
\usepackage{graphicx}
\usepackage{epsfig}
\newtheorem{lemme}{Lemma} 
\newtheorem{theorem}{Theorem}

\newtheorem{proposition}{Proposition} 
\def\R{{\mathbb R}} 
 
\def\Z{{\mathbb Z}} 

\def\N{{\mathbb N}}

\def\V{\text{Vol}}

\def\s{\text{sys}}

\begin{document}

\title{Volume entropy, systole and stable norm on graphs}

\author{Florent Balacheff \footnote{F. Balacheff  UMR 5149, Institut de Math\'ematiques et de Mod\'elisation de Montpellier - Universit\'e de Montpellier Case Courrier 051
-    place Eug\`ene Bataillon 34095 Montpellier Cedex 5, France (email : balachef@math.univ-montp2.fr)}}

\maketitle

\begin{abstract}
We study some new isoperimetric inequalities on graphs. We etablish a relation between the volume entropy (or asymptotic volume), the systole and the first Betti number of weighted graphs. We also find bounds for the volume, associated to some special measure, of the unit ball for the stable norm of graphs.
\end{abstract}

\bigskip
\noindent {\it Mathematics Subject Classification (2000)} : 05C35, 37B10 .

\section{Introduction}

Some isoperimetric inequalities on manifolds can sometimes be extend to graphs. For example, the asymptotic behaviour of the systolic volume of a surface in terms of its genus, investigated in \cite{bussar} and \cite{grom1}, has a $1$-dimensional analog valid on weighted graphs \cite{bolsze}.  Namely, let $(\Gamma,w)$ be a weighted graph (see subsection 2.1 for a precise definition). The {\it volume} of $(\Gamma,w)$ (or {\it size}), denoted by $\V(\Gamma,w)$, is the sum of the weight of its edges 
$$
\V(\Gamma,w)=\sum_{e \in E} w(e).
$$
 The {\it systole} of $(\Gamma,w)$ (or {\it weighted girth}) is defined as 
$$
\s(\Gamma,w)=\inf \{l_w(\gamma) \mid \gamma \text{ non trivial cycle of } \Gamma \},
$$
 where the length of a cycle $\gamma$, noted $l_w(\gamma)$, is the sum of the weights of its edges.

With this two quantities, we can define the {\it systolic volume} of $\Gamma$ as 
$$
\sigma(\Gamma)=\inf_{w} \frac{\V(\Gamma,w)}{\s(\Gamma,w)}. 
$$
where the infimum is taken over all the weight functions on the graph $\Gamma$. 

\bigskip
\noindent {\bf Theorem (Bollob\'as \& Szemer\'edi \cite{bolsze})} {\it  Let $\Gamma$ be a finite graph with first Betti number $b \geq 3$. We have :
$$
\sigma(\Gamma) \geq \frac{3\ln 2}{2} \frac{b-1}{\ln (b-1)+\ln \ln (b-1)+4\ln 2 -\ln \ln 2}. 
$$
}

\bigskip

This lower bound gives an asymptotic lower estimate in terms of the first Betti number $b$:
$$
\sigma(\Gamma) \geq \frac{3 \ln 2}{2} \frac{b}{\ln b} +o\left(\frac{b}{\ln b}\right ). \eqno (1.1)
$$

On the other hand, for each $b \geq 2$, we can construct a weighted graph $(\Gamma,w)$  (see \cite{babbal}) of first Betti number  $b$ such that
$$
\sigma(\Gamma,w) \leq 8 \ln 2 \frac{b}{\ln b},
$$
where $\sigma(\Gamma,w)=\V(\Gamma,w) / \s(\Gamma,w)$.

Thus, the estimate (1.1) gives a good asymptotic behaviour of the systolic volume in terms of the first Betti number.  This is the analog of an estimate on surfaces, which states that the systolic volume of a surface $S$ of genus $g$, denoted by $\sigma(S)$ (see \cite{grom1} for a precise definition), has the following lower bound
$$
\sigma(S) \geq C \frac{g}{(\ln g)^2}, \eqno (1.2)
$$
where $C$ is some positive constant.
In \cite{bussar}, the authors construct a metric on surfaces of genus $g$ for which the systolic volume is asymptotically close to (1.2).

\bigskip

 In this article, we are interested in the normalization of the volume entropy by the systole. Let $(\tilde{\Gamma},\tilde{w})$ be the universal (weighted) covering of $(\Gamma,w)$. Fix $x_0 \in \Gamma$ and $\tilde{x}_0 \in \tilde{\Gamma}$ a lift of $x_0$. The {\it volume entropy} (or {\it asymptotic volume}) of $(\Gamma,w)$ is defined as
$$
h_{vol}(\Gamma,w)=\lim_{R \rightarrow +\infty} \frac{\ln(\V_{\tilde{w}}B(\tilde{x}_0,R))}{R} \eqno (1.3)
$$
where $\V_{\tilde{w}} B(\tilde{x}_0,R)$ is the volume of the ball centered at $\tilde{x}_0$ with radius $R$ in $(\tilde{\Gamma},\tilde{w})$. Since the weighted graph $(\Gamma,w)$ is compact,  the limit in (1.3) exists and does not depend on the point $x_0 \in \Gamma$ and its lift (see \cite{man}). 

The product $h_{vol}(\Gamma,w).\s(\Gamma,w)$ is invariant under scaling and has been studied in \cite{sab} for surfaces. The author has proved that this quantity is bounded from above for each surface. We will prove  an upper bound on this product for weighted graphs, which is asymptotically optimal.

\bigskip
\noindent {\bf Theorem 1} {\it Let $(\Gamma,w)$ be a weighted graph with first Betti number equal to $b$. Then
$$
h_{vol}(\Gamma,w).\s(\Gamma,w) \leq 2 \ln (8b^3-1).
$$ 
}

We state an analogue result for topological Markov chains associated to a $n \times n$ matrix $A$, relating the topological entropy $h_{top}$ and the smallest period of a periodic orbit $T_{min}$.

\bigskip

\noindent {\bf Proposition} {\it For each topological Markov chain $(\Sigma_A,\phi_A)$, 
$$
h_{top}(\Sigma_A,\phi_A).T_{min}(\Sigma_A,\phi_A)\leq \ln b_A,
$$
where $b_A=\sum_{i,j=1}^n A_{ij}-n +1$.
}
\medskip

This inequality is easy to prove. We present here because, as far as the author is aware, no reference of this statement exists in the literature.

\bigskip

In an another direction, we can associate to a weighted graph an object called the {\it stable norm} on the real first homology space (see subsection 3.1 for a definition). We can define a natural measure $\mu_w$ on $H_1(\Gamma,\R)$, which allows us to define the volume of the {\it stable ball} ${\cal B}_{st}(\Gamma,w)$, defined as the unit ball of the stable norm. For usual combinatorial graphs, which can of course be identified with a weighted graph in which all the edges have weight $1$ (in this case, we forget the indice $w$), we obtain the following inequalities.

\bigskip

\noindent {\bf Theorem 2} {\it Let $\Gamma$ be a graph with first Betti number $b$. Then
$$
\frac{2^b} {b!}\left( \frac{b}{\V(\Gamma)}\right)^b \leq \mu ({\cal B}_{st}(\Gamma)) \leq \frac{2^b}{b!}, 
$$
The two cases of equality are attained by the bouquet of $b$ circles $\bigvee_{i=1}^b S^1_i$.
}
\bigskip

We also obtain an inequality on weighted graphs relating the volume of the graph to the volume of the stable ball.

\bigskip

\noindent {\bf Theorem 3} {\it For every weighted graph $(\Gamma,w)$ with first Betti number $b$,
$$
\mu_w({\cal B}_{st}(\Gamma,w)).\V(\Gamma,w)^{b/2} \geq \omega_b,
$$
where $\omega_b$ is the volume of the Euclidean unit ball of $\R^b$.
}
\bigskip

The first part of this paper is dedicated to the proof of theorem {\bf 1} and the proposition above. Some related problems are also considered. In the second part, we recall some definitions and show theorems {\bf 2} and {\bf 3}.

\section{Volume entropy, systole and scale of graphs}

\subsection{Preliminaries}

By a finite graph, we mean a finite non-oriented multigraph (we allow multiple edges and loops). For a finite graph  $\Gamma$, we denote by $V$ the set of vertices and by $E$ the set of edges. Each element of $E$ is an element of $V \times V$. A {\it weighted graph} is a pair $(\Gamma,w)$, where $\Gamma$ is a finite graph and $w : E \rightarrow \mathbb R^+$ is a weight function. For $e \in E$, we call $w(e)$ the weight of $e$. From now on, all graphs considered are finite.

We recall the definition of the exponential growth rate of a group (see \cite{harp}), that we will use in the proof of theorem {\bf 1}, and in the proof of proposition {\bf 3} and {\bf 5}.

Let $G$ be a group of finite presentation and $\Sigma$ be a finite generating set. We define the {\it algebraic length} of an element $\alpha$ of $G$ with respect to $\Sigma$ as the smallest integer $k$ such that $\alpha=\alpha_1\ldots\alpha_k$, where $\alpha_i \in \Sigma \cup \Sigma^{-1}$. It is denoted by $|\alpha|_\Sigma$.

The {\it exponential growth rate} of $G$ with respect to the system $\Sigma$ is defined as
$$
\omega(G,\Sigma)=\lim_{R\rightarrow +\infty} \sqrt[k]{N_\Sigma(R)} 
$$
where $N_\Sigma(R)=\text{card}\{\alpha \in G \mid |\alpha|_\Sigma\leq R\}$ is the cardinal of the ball of radius  $R$ of $(G,|.|_\Sigma)$ centered at its origin. For a group of finite presentation $G$ and a finite generating set $\Sigma$ of $G$,
$$
\omega(G,\Sigma)\leq 2.\text{card}(\Sigma)-1. \eqno (1.3)
$$

\subsection{Volume entropy and systole of regular graphs}

In this subsection, we etablish an upper bound on the normalization of the volume entropy by the systole valid on regular graphs which is better than the bound of theorem {\bf 1}. Recall that the {\it valence} of a vertex of a graph is the number of incident edges at this vertex. Let $v\geq 2$ be an integer. We say that a graph $\Gamma$ is {\it regular} of valence $v$ if the valence of each vertex is constant equal to $v$. Recall that a (usual combinatorial) graph is naturally identified with a weighted graph in which all the edges have weight $1$. Remark that in this case  $\V(\Gamma)=\text{card}(E)$.

We will show the following result.

\begin{proposition}
Let $\Gamma$ be a regular graph with first Betti number $b$. Then
$$
h_{vol}(\Gamma).\s(\Gamma) \leq 3 \ln b. \eqno (1.4)
$$
\end{proposition}

\noindent {\bf Proof of the proposition.} Let  $v$ be the valence of $\Gamma$. We have 
$$
h_{vol}(\Gamma)=\ln(v-1). \eqno (1.5)
$$
To see this, fix some vertex $\tilde{x}_0$ in the universal covering $\tilde{\Gamma}$. As $\tilde{\Gamma}$ is an infinite regular tree of valence $v$, 
\begin{eqnarray}
\nonumber \V(B(\tilde{x}_0,R))& = & v(1+(v-1)+\ldots +(v-1)^{R-1})\\
\nonumber&=&v \frac{(v-1)^R-1}{v-2}\\
\nonumber
\end{eqnarray}
for each positive integer $R$. We deduce (1.5).

\bigskip

In the case $b=1$, we obtain $v=2$ and so $h_{vol}(\Gamma)=0$. The inequality (1.4) is then trivial.

In the case $\s(\Gamma)=1$, we have $b=\text{ card}(E)-\text{ card}(V)+1$, and as, by elementary considerations, $2.\text{card}(E)=v.\text{card}(V)$, and  $\text{card}(V) \geq 1$, we get $b \geq v/2$. Therefore,
\begin{eqnarray}
\nonumber h_{vol}(\Gamma)&=&\ln(v-1)\\
\nonumber &\leq&\ln(2b-1)\\
\nonumber &\leq& 3 \ln b,\\
\nonumber
\end{eqnarray}
and (1.4) follows in this case.

Now suppose that $b>1$ and $\s(\Gamma)>1$. We will show the following lemma.

\begin{lemme}
$$
 \s(\Gamma) \leq \frac{3 \ln b}{\ln(v-1)}. \eqno(1.6)
$$
\end{lemme}

\noindent {\bf Proof.} For all $R < \s/2$, the ball centered at any point $x$ in $\Gamma$ of radius $R$ is a tree. Thus, the calcul of the volume of a ball centered in a vertex of radius $[\s(\Gamma)/2]$ gives the following estimate (compare with \cite{bol}, p.14)
$$
\text{card}(E) \geq v\frac{(v-1)^{[\s(\Gamma)/2]}-1}{v-2}. 
$$
With $\text{card}(E)=v(b-1)/(v-2),$ we deduce   
$$
[\s(\Gamma)/2] \leq \frac{\ln b}{\ln (v-1)}.
$$

Then 
$$
\s(\Gamma) \leq  1 + 2\frac{\ln b}{\ln (v-1)} \leq 3\frac{\ln b}{\ln (v-1)},
$$
and we are done.

\bigskip

\noindent Now we combine inequalities (1.5) and (1.6) to get the inequality (1.4). \hspace{\stretch{1}} $\Box$\\

\bigskip

The asymptotic behaviour of (1.4) is optimal when $b$ goes to infinity. It is realized by the bouquet of $b$ circles viewed as the graph with one vertex and $b$ loops of weight $1$.

\bigskip

In the case of a graph provided with some control on the valence of its vertices, we easily obtain a lower and an upper bound on the volume entropy. We denote by $v(s)$ the valence of a vertex $s \in V$.

\begin{proposition}
Let $\Gamma=(V,E)$ be a graph. Suppose that there exists two integers $2 \leq \delta \leq \Delta$ such that for all $s \in V$, $\delta \leq v(s)\leq \Delta$. Then
$$
\ln(\delta -1) \leq h_{vol}(\Gamma) \leq \ln( \Delta-1).
$$
\end{proposition}

\subsection{Volume entropy and systole of weighted graphs}

We now generalize (1.4) to weighted graphs.

\begin{theorem}
Let $(\Gamma,w)$ be a weighted graph with first Betti number equal to $b$. Then
$$
h_{vol}(\Gamma,w).\s(\Gamma,w) \leq 2 \ln (8b^3-1). \eqno (1.7)
$$
\end{theorem}

\noindent {\bf Proof.} The proof involves techniques coming from \cite{svar}. Fix a fundamental domain $D$ of $\tilde{\Gamma}$ and a point $\tilde{x}_0$ in $D$ such that points of the boundary are not vertices  of $\tilde{\Gamma}$. $D$ is a tree and denote by $\{\tilde{y}_1,\ldots,\tilde{y}_m\}$ the boundary. We have $m\leq 2b$. Denote by $p : \tilde{\Gamma} \rightarrow \Gamma$ the universal covering projection  and by $\{z_1,\ldots,z_k\}$ the image of $\{\tilde{y}_1,\ldots,\tilde{y}_m\}$ under $p$. It is obvious that $k\leq b$. Denote by $x_0$ the projection of $\tilde{x}_0$.

Let $s=\s(\Gamma,w)/2$. Denote by  $p(s)$ the minimal number of translated domains $\gamma.D$ under the action of $\pi_1(\Gamma,x_0)$ needed to form a neighbourhood of $D$ such that each point not belonging to this neighbourhood is at a distance more than $s$ from $D$. We will find an upper bound of $p(s)$. Fix a vertex $\tilde {y}_i$ in the boundary  and enumerate the paths starting at $\tilde{y}_i$ of length $s$ going outside $D$. The number of these paths is less than the number of paths starting at $p(\tilde{y}_i)$ of length $s$, which is less than $2b$. It is clear that each of these paths passes at most once through each element of $\{z_1,\ldots,z_k\}$ (as $s=\s(\Gamma,w)/2$). If we consider a path $\tilde{c}$ of $\tilde{\Gamma}$ starting at $\tilde{y}_i$ of length $s$, the number of translated domain that it passes through is exactly the number of points in $\{z_1,\ldots,z_k\}$  that belongs to $p(\tilde{c})$. Thus, this number is bounded from above by $k\leq b$. The number of translated domains of non null intersection with a path starting from $\tilde{y}_i$ with length $s$ is then bounded from above by $2b^2$.

 As there are $m$ elements in the boundary of $D$, and $m\leq 2b$, we get  $p(s) \leq 4b^3$.

\bigskip
Denote by $\Sigma(s)=\{\gamma_i\}_{i=1}^{p(s)}$ a generating set of $\pi_1(\Gamma,x_0)$ such that for each $x \notin \cup_{i=0}^{p(s)} \gamma_i.D$ we have $d_{\tilde{w}}(x,D)\geq s$ (here $\gamma_0$ is by convention the neutral element).

We will estimate the volume of the ball centered at $\tilde{x}_0 \in D$ with radius $ns$. We can easily show that
$$
B^{\tilde{\Gamma}}(x_0,ns) \subset B(\Sigma(s),n),
$$
where $ B(\Sigma(s),n) :=  \cup \{ \gamma. D \mid |\gamma|_{\Sigma(s)} \leq n \}$.

We deduce 
$$
\V_{\tilde{w}}(B^{\tilde{\Gamma}}(\tilde{x}_0,ns)) \leq  \V_{\tilde{w}}(\Gamma).N_{\Sigma(s)}(n),
$$
and so, by (1.3),
$$
s.h_{vol}(\Gamma,w)\leq \ln (2.p(s)-1).
$$
Inequality (1.7) follows. \hspace{\stretch{1}} $\Box$\\

\bigskip

Inequality (1.7) yields the asymptotic behaviour
$$
h_{vol}(\Gamma,w).\s(\Gamma,w) \lesssim 6\ln b. 
$$

We say that a weighted graph $(\Gamma,w)$ has a {\it systolic basis} if there exist $x_0 \in V(\Gamma)$ and  a generating set $\Sigma_0=\{\gamma_i\}_{i=1}^b$ of $\pi_1(\Gamma,x_0)$ such that $l_w(\gamma_i)=\s(\Gamma,w)$ for $i=1,\ldots,b$. For weighted graph with systolic basis, we have the following lower bound.

\begin{proposition}
Let $(\Gamma,w)$ be a weighted graph with first Betti number $b$. Suppose that $(\Gamma,w)$ has a systolic basis. Then
$$
h_{vol}(\Gamma,w).\s(\Gamma,w) \geq \ln(2b-1), \eqno (1.8)
$$
and the equality case is attained by the bouquet of $b$ circles $\bigvee_{i=1}^b S^1_i$.
\end{proposition}

\noindent {\bf Proof.} Let $x_0 \in \Gamma$ be the vertex of $\Gamma$ which is the base-point of the systolic basis. We fix $\tilde{x}_0 \in \tilde{\Gamma}$ a lift of $x_0$. The geometric length of an element $\alpha \in \pi_1(\Gamma,x_0)$ is defined by
\begin{eqnarray}
\nonumber |\alpha|_w & = & \min \{l_w(\gamma) \mid \gamma \text{ cycle based at } x_0 \text{ representing } \alpha \}\\
\nonumber &=&d_{\tilde{w}}(\tilde{x}_0,\alpha.{\tilde{x}_0}).\\
\nonumber
\end{eqnarray} 
This geometric length does not depend on the lift $\tilde{x}_0$.
\bigskip

We then have (see \cite{sab})
$$
h_{vol}(\Gamma,w)=\lim_{R \rightarrow +\infty} \frac{\ln (N_w(R))}{R}
$$
where $N_w(R)=\text{card} \{\alpha \in \pi_1(\Gamma,x_0) \mid |\alpha|_w \leq R\}$ is the cardinal of the ball of $(\pi_1(\Gamma,x_0),|.|_w)$ with radius $R$ centered at the origin.

We have easily
\begin{lemme}
Let $(\Gamma,w)$ be a weighted graph and $\Sigma$ be a generating set of $G=\pi_1(\Gamma,x_0)$. If there exist $\lambda, \mu >0$ such that
$$
\lambda.|.|_\Sigma \leq |.|_w \leq \mu.|.|_\Sigma,
$$
then
$$
\frac{1}{\mu} \ln \omega(G,\Sigma) \leq h_{vol}(\Gamma,w) \leq \frac{1}{\lambda} \ln \omega(G,\Sigma). \eqno (1.9)
$$
\end{lemme}

Now, as $(\Gamma,w)$ has a systolic basis $\Sigma_0$ at $x_0$, for every $\alpha \in \pi_1(\Gamma,x_0)$, 
$$
|\alpha|_w \leq |\alpha|_{\Sigma_0}.\s(\Gamma,w).
$$
We immediately deduce (1.8).

The case of equality is attained by the bouquet $\vee_{i=1}^b S^1_i$. \hspace{\stretch{1}} $\Box$\\

\bigskip

\subsection{Volume entropy and scale of graphs}

Let $(\Gamma,w)$ be a weighted graph. A {\it chain} of $\Gamma$ is a path such that the valence of each intermediate vertex is $2$.

We define the {\it microscopic scale} of $(\Gamma,w)$ as 
$$
C_{min}(\Gamma,w)=\min \{l_w(C) \mid C \text{ chain of } \Gamma \},
$$
and the {\it macroscopic scale} of $(\Gamma,w)$ as 
$$
C_{max}(\Gamma,w)=\max \{l_w(C) \mid C \text{ chain of } \Gamma \}.
$$
The aim of this subection is to prove isoperimetric inequalities involving the volume entropy and the scale of weighted graphs. 

\begin{proposition}
Let $(\Gamma,w)$ be a weighted graph such that $v(s)\leq 3$ for every $s \in V$.

Then
$$
\frac{\ln 2}{C_{max}(\Gamma,w)} \leq h_{vol}(\Gamma,w)\leq \frac{\ln 2}{C_{min}(\Gamma,w)}. \eqno (1.10)
$$
\end{proposition}

\noindent {\bf Proof.} Let $T_3$ be the infinite regular tree of valence $3$. Fix a vertex $v$ in $T_3$ and $v'$ in $\Gamma$. We denote by $w_{min}$ (respectively $w_{max}$) the constant weight function defined on $T_3$ equal to $C_{min}(\Gamma,w)$ (respectively $C_{max}(\Gamma,w)$). Then, for every $R>0$,
 $$
\V_{w_{max}} (B^{T_3}_{w_{max}}(v,R)) \leq \V_w(B^{\tilde{\Gamma}}(v',R))\leq \V_{w_{min}}(B^{T_3}_{w_{min}} (v,R))
$$
and so
$$
3 (2^{[R/C_{max}]}-1).C_{max} \leq \V_w(B^{\tilde{\Gamma}}(v',R))\leq 3 (2^{[R/C_{min}]+1}-1).C_{min}.
$$
We deduce (1.10). \hspace{\stretch{1}} $\Box$\\

We can prove a stronger result for the normalization of the volume entropy by the minimum scale.
\begin{proposition}
Let $(\Gamma,w)$ be a weighted graph of first Betti number $b$. Then
$$
h_{vol}(\Gamma,w).C_{min}(\Gamma,w)\leq \ln (2b-1). \eqno (1.11)
$$
The equality case is attained by the bouquet of $b$ circles $\bigvee_{i=1}^b S^1_i$.
\end{proposition}

\noindent {\bf Proof.} Fix a vertex $x_0$ and let $\Sigma=\{\gamma_1,\ldots,\gamma_b\}$ be a minimal generating set of $\pi_1(\Gamma,x_0)$. Choose a maximal tree $T$ of $\Gamma$ containing $x_0$, and denote by $\{e_i\}_{i=1}^b$ the edges of $\Gamma \setminus T$ such that  $e_i \in \gamma_i$ for $i=1,\ldots,b$ and $e_j \notin \gamma_i$ for $j\neq i$. Also denote by $p : \Sigma \cup \Sigma^{-1} \rightarrow \{e_1,\ldots,e_b\}$ the map defined by $p(\gamma_i)=p(\gamma_i^{-1})=e_i$. For every $\alpha \in \pi_1(\Gamma,x_0)$, we choose a reduced form  $\alpha=\alpha_1\ldots\alpha_k$ where $\alpha_j \in \Sigma \cup \Sigma^{-1}$ and $k=|\gamma|_\Sigma$. For every cycle $\beta$ homotopic to $\alpha$, $\beta$ may be written as the concatenation of paths
$$
[\beta_1,p(\alpha_1),\beta_2,p(\alpha_2),\ldots,p(\alpha_k),\beta_{k+1}]
$$
where $\beta_i$ are paths of $\Gamma$.
As $l_w(e_i) \geq C_{min}(\Gamma,w)$ for $i=1,\ldots,k$, we get
$$
l_w(\beta) \geq |\alpha|_\Sigma.C_{min}(\Gamma,w),
$$
so
$$
|\alpha|_w \geq |\alpha|_\Sigma.C_{min}(\Gamma,w).
$$
From inequality (1.9), we deduce (1.11). \hspace{\stretch{1}} $\Box$\\

\subsection{Entropy and topological Markov chain}

We show here the proposition {\bf 6} stated in the introduction. We recall first some definitions (see \cite{shub}). Let $n$ be a positive integer. We denote by $\N_n$ the set $\{1,\ldots,n\}$ and $\Sigma(n)$ the product space $\N_n^\Z$. The product topology is then induced by the following metric on $\Sigma(n)$
$$
d(a,b)=\sum_{i=-\infty}^{+\infty} \frac{d_i(a,b)}{2^{2|i|+1}},
$$
where $d_i(a,b)$ is equal to $0$ if $a_i=b_i$, or to $1$ otherwise. Note that the sequence $(a^j)_{j \in \N}$ converges if and only if for all $i \in \Z$ the sequence $(a_i^j)_{j \in \N}$ converges.
Let $\phi$ be the homeomorphism of $\Sigma(n)$ defined as $(\phi(a))_i=a_{i+1}$. The homeomorphism $\phi$ is called the {\it shift}.

Denote by $M_n$ the space of the $n\times n$ matrices all of whose entries are $0$ or $1$. If $A \in M_n$, we define
$$
\Sigma_A=\{ a \in \Sigma(n) \mid A_{a_i,a_{i+1}}=1\}.
$$
The set $\Sigma_A$ is a closed $\phi$-invariant subspace of $\Sigma(n)$. We denote by $\phi_A : \Sigma_A \rightarrow \Sigma_A$ the restriction of $\phi$ to $\Sigma_A$. The pair $(\Sigma_A,\phi_A)$ is called a {\it subshift of finite type} or {\it topological Markov chain}.

The {\it topological entropy} of a dynamical system (see \cite{adkomc} for a definition) is denoted by $h_{top}$. In the case of a topological Markov chain $(\Sigma_A,\phi_A)$, we deduce from \cite{adkomc} that 
$$
h_{top}(\Sigma_A,\phi_A)= \lim_{k \rightarrow +\infty} \frac{\ln (N(\Sigma_A,k))}{k},
$$
where $N(\Sigma_A,k)$ is the number of open sets in the following cover
$$
{\cal A}_k =\{\{a \in \Sigma(n) \mid   \forall j=0, \ldots,k-1 ,  a_j=i_j\} \mid  \{i_j\} \in \N^k_n \}.
$$
We introduce also $N_i(\Sigma_A,k)$ the cardinal of open sets in ${\cal A}_k\cap \{a \mid a_0=i\}$. We have
$$
N(\Sigma_A,k)=\sum_{i=0}^k N_i(\Sigma_A,k).
$$

It is a classic result that the topological entropy of $(\Sigma_A,\phi_A)$ is exactly the natural logarithm of the spectral radius of $A$.

\bigskip

From a different point of view, $\Sigma_A$ can be identified with the space of bi-infinite paths of an oriented graph $\Gamma_A$ whose vertices lie in $\N_n$ and edges are the pairs $(i,j)$ of vertices such that $A_{i,j}=1$. We define the {\it minimal period} of $(\Sigma_A,\phi_A)$, denoted  by $T_{min}(\Sigma_A,\phi_A)$, as the smallest period of a periodic point of the dynamical system $(\Sigma_A,\phi_A)$. It coincides with the shortest length of an oriented cycle of $\Gamma_A$. We denote by $b_A$ the first Betti number of $\Gamma_A$ which can be expressed as $b_A=\sum_{i,j=1}^n A_{ij}-n+1$. 

\begin{proposition}
For each topological Markov chain $(\Sigma_A,\phi_A)$, 
$$
h_{top}(\Sigma_A,\phi_A).T_{min}(\Sigma_A,\phi_A)\leq \ln b_A. \eqno (1.12)
$$
The equality case is attained for each $b \geq 1$ by the following Markov chain defined as the  matrix of size $b+1$
$$
\left\{ \begin{array}{l}
A_{1,j}=A_{j,1}=1 \text{ for } j \neq 1,\\
A_{i,j}=0 \text{ otherwise}.\\
\end{array}
\right. \eqno (1.13)
$$
\end{proposition}

\noindent {\bf Proof.} 
Fix $t=T_ {min}(\Sigma_A,\phi_A)$. For each vertex $v \in \N_n$ of $\gamma_A$, the number of oriented paths starting at $v$ of length $t$ is at most $b_A$.

We obtain  
$$
N_i(\Sigma_A,nt) \leq b_A^n,
$$
and so 
$$
N(\Sigma_A,nt) \leq n.b_A^n.
$$
Then
$$
t.h(\Sigma_A,\phi_A) \leq \ln b_A.
$$
and we get the claim (1.12). It is easy to compute that the topological Markov chain defined by (1.13) realize the equality case. \hspace{\stretch{1}} $\Box$\\

\section{Stable norm and volume of graphs}

\subsection{Unit ball of the stable norm}

Let $\Gamma=(V,E)$ be a finite graph with first Betti number $b$. For the following classical notions, we refer to \cite{span}. The graph $\Gamma$ is a simplicial complex of dimension $1$, whose $1$-simplices are the edges and $0$-simplices are the vertices. We denote by $C(\Gamma,\mathbb R)$ the {\it oriented chain complex} of the simplicial complex $\Gamma$ which is, in this 1-dimensional context and after a choice of orientation for each edge, the real vectorspace $\oplus_{e \in E} \mathbb R e$ of dimension $|E|$. Since the homology groups of dimension at least two are null, the  1-dimensional homology group with real coefficients $H_1(\Gamma,\R)$ is embedded in $C(\Gamma,\mathbb R)$ as a subspace of dimension $b$. Finally, the 1-homology group with integer coefficients $H_1(\Gamma,\Z)$ is identified with a lattice of the subspace $H_1(\Gamma,\R)$ (compare with \cite{bacharnag}).

Let $w$ be a weight function on $\Gamma$. Denote by the sequence $\{e_i\}_{i=1}^{|E|}$ the edges of $E$ and by $w_i=w(e_i)$ the weight of each edge.  We define the following scalar product denoted by  $<.,.>_w$
$$
<e_i,e_j>_w=w_i \delta_{ij},
$$ 
for $i,j=1,\ldots,|E|$, where $\delta_{ij}$ is the Kronecker symbol.

 Recall now the definition of the stable norm defined on $H_1(\Gamma,\R)$ (see \cite{grom}). For an element $u \in H_1(\Gamma,\R)$, 
$$
|u|_{st}=\lim_{n \rightarrow +\infty} \frac{1}{n} |\overline{nu}|_w, \eqno (2.1)
$$
where $\overline{nu}$ is the point of $H_1(\Gamma,\Z)$ closest (for an arbitrary fixed Euclidean metric) to $nu$. Recall that $| \cdot |_w$ is the geometric length on $\pi_1(\Gamma,x_0)$. The limit in (2.1) always exists and is independent on $x_0$.

We define the {\it unit ball of the stable norm} of $(\Gamma,w)$ by
$$
{\cal B}_{st}(\Gamma,w)=\{u \in H_1(\Gamma,\R) \mid |u|_{st}\leq 1\}.
$$

\medskip

Denoted by $\mu_w$ the Haar measure on $H_1(\Gamma,\R)$ induced by the restriction of the scalar product $<.,.>_w$ to this subspace. When we consider a combinatorial graph, we forget the indice $w$ in the notations $<.,.>_w$ and $\mu_w$. For a combinatorial graph $\Gamma$, the measure $\mu$ coincides with the restriction to $H_1(\Gamma,\R)$ of the Lebesgue measure of $\R^{|E|}$.

\subsection{Inequalities for combinatorial graphs}

\begin{theorem}
Let $\Gamma$ be a graph with first Betti number $b$. We have the following inequalities
$$
\frac{2^b}{b!}\left( \frac{b}{\V(\Gamma)}\right)^b \leq\mu({\cal B}_{st}(\Gamma)) \leq \frac{2^b}{b!}. \eqno (2.2)
$$
The two cases of equality are attained by the bouquet of $b$ circles $\bigvee_{i=1}^b S^1_i$.
\end{theorem}

\noindent {\bf Proof.} We have
$$
{\cal B}_{st}(\Gamma)=B_1^k \cap H_1(\Gamma,\R), \eqno (2.3)
$$
where $k=|E|$ and $B_1^k=\{(x_i) \in \R^k \mid \sum_{i=1}^k |x_i|\leq 1\}.$

Indeed, we can show that the stable norm agrees with the norm $|.|_1$ defined by $|u|_1=\sum_{i=1}^k |u_i|$ for all $u=\sum_{i=1}^k u_i e_i  \in H_1(\Gamma,\R)$. Namely,
\begin{eqnarray}
\nonumber |u|_{st}&=&\lim_{n \rightarrow +\infty} \frac{1}{n} |\overline{nu}|_w\\
\nonumber&=&\lim_{n \rightarrow +\infty} \frac{1}{n} |\overline{nu}|_1\\
\nonumber&=&|u|_1.\\
\nonumber
\end{eqnarray}

\medskip

We can find in  \cite{meypaj} the following estimate. For every  $b$-plane $P^b$  in $\R^k$
$$
\mu_k(B_1^k \cap P^b) \geq \mu_k(B_1^k)^{b/k},
$$
where $\mu_n$ is the canonical volume of $\R^n$.
We deduce with (2.3) that
$$
\mu(B_{st}(\Gamma)) \geq \mu_k(B_1^k)^{b/k}.
$$
Since  $\mu_n(B_1^{n})=2^n/n!$ for every positive integer $n$, we have 
$$
\mu(B_{st}(\Gamma)).k^b \geq \frac{2^b.k^b}{(k!)^{b/k}}.
$$
With $k\geq b$, we get
$$
\mu(B_{st}(\Gamma')).k^b \geq \frac{(2b)^b}{b!}.
$$
As $k=|E|=\V(\Gamma)$, the left inequality of (2.2) is then proved.

For the upper bound, we start with an other estimate obtained in \cite{meypaj}. For every $b$-plane $P^b$  in $\R^k$,
$$
\mu_k(B_1^k \cap P^b) \leq \mu_b(B_1^b).
$$
 So
$$
\mu(B_{st}(\Gamma)) \leq \mu_b(B_1^b),
$$
which gives the right inequality of (2.2).

\bigskip

For the bouquet of $b$ circles (viewed as the combinatorial graph with one vertex and $b$ loops)
$$
\mu(B_{st}(\bigvee_{i=1}^b S^1_i))=\frac{2^b}{b!}
$$
and $\V(\bigvee_{i=1}^b S^1_i)=b$, so $\bigvee_{i=1}^b S^1_i$ realize the equality case.

\hspace{\stretch{1}} $\Box$\\

\bigskip
\noindent {\bf Remark}
For a regular graph $\Gamma$ of valence $v\geq 3$, we have 
$$
\left( \frac{v-2}{v} \right)^b \frac{2^b}{b!} \leq \mu_w({\cal B}_{st}(\Gamma)) \leq \frac{2^b}{b!}.
$$

\bigskip

\subsection{Estimate for weighted graphs}

\begin{theorem}
For every weighted graph $(\Gamma,w)$ with first Betti number $b$,
$$
\mu_w({\cal B}_{st}(\Gamma,w)).\V(\Gamma,w)^{b/2} \geq \omega_b, \eqno (2.4)
$$
where $\omega_b$ is the volume of the Euclidean unit ball of $\R^b$.
\end{theorem}

\noindent {\bf Proof.} Suppose that the result holds for graphs. Let $(\Gamma,w)$ be a weighted graph with first Betti number $b$.  For each $\epsilon >0$, we can find a weight function $w_\epsilon$ close enough to $w$ in the sense of the $C^0$ topology such that  
$$
|\mu_{w_\epsilon}({\cal B}_{st}(\Gamma,w_\epsilon))-\mu_w({\cal B}_{st}(\Gamma,w))|<\epsilon,
$$
$$
|\V(\Gamma,w_\epsilon)-\V(\Gamma,w)|<\epsilon,
$$
and such that  $w_\epsilon(e)$ is rational for every $e \in E$. Fix an integer $\lambda$ such that $\lambda w_\epsilon(e) \in \N^\ast$ for all $e \in E$. We have 
$$
\mu_{\lambda  w_\epsilon}({\cal B}_{st}(\Gamma,\lambda  w_\epsilon)).\V(\Gamma,\lambda w_\epsilon)^{b/2}=\mu_{w_\epsilon}({\cal B}_{st}(\Gamma,w_\epsilon)).\V(\Gamma,w_\epsilon)^{b/2}.
$$ 
Choose an enumeration $\{e_i\}_{i=1}^{|E|}$ of the edges of $\Gamma$. If we subdivise each edge $e_i$ of $\Gamma$ in $k_\epsilon(i)=\lambda w_\epsilon(e_i)$ edges denoted by $e'_{i,j}$ for $j=1,\ldots,k_\epsilon(i)$, we  get a graph $\Gamma'_\epsilon$ with first Betti number $b$, which is isometric to $(\Gamma,\lambda w_\epsilon)$ when it is endowed with the trivial weight function $w=1$. Denote by $f$ this isometry and by $E'$ the set of the edges of $\Gamma'_\epsilon$. We have
$$
\mu({\cal B}_{st}(\Gamma_\epsilon')) = \mu_{\lambda  w_\epsilon}({\cal B}_{st}(\Gamma,\lambda w_\epsilon)). \eqno (2.5)
$$
To see this, observe that the isometry $f$ induces a linear homomorphism $F$ from $C(\Gamma,\R)$ to $C(\Gamma_\epsilon',\R)$ (and an isomorphism between  $H_1(\Gamma,\R)$ and $H_1(\Gamma'_\epsilon,\R)$) which satisfies
$$
F(e_i)=\sum_{j=1}^{k_\epsilon(i)} e'_{i,j},
$$
for $i=1,\ldots,|E|$.
The familly 
$$
\left \{\frac{1}{\sqrt{k_\epsilon(i)}} e_i\right \}_{i=1}^{|E|}
$$
is an orthonormal basis for the scalar product $<.,.>_{\lambda w_\epsilon}$ of $\R^{|E|}$ and 
$$
\left \{\frac{1}{\sqrt{k_\epsilon(i)}} \sum_{j=1}^{k_\epsilon(i)} e'_{i,j} \right \}_{i=1}^{|E|}.
$$
an orthonormal basis of $I(C(\Gamma,\R))$ for $<.,.>$ in $\R^{|E'|}$. We then find that the map $F_{|H_1(\Gamma,\R)}$ expressed in these orthonormal basis is the identity map, so we get the claim (2.5).

\medskip

This construction can be realized for every $\epsilon>0$. So if the result holds for graphs, it holds for every weighted graph.

\medskip

Let $\Gamma'=(E',V')$ be a graph with first Betti number $b$ and denote $k=|E'|$. For every $u \in \R^k$, 
$$
|u|_1 \leq k.|u|^2_2
$$ 
where $|u|_2=\sqrt{\sum_{i=1}^{k} u_i^2}$ is the Euclidean norm. So 
$$
\mu({\cal B}_{st}(\Gamma')) \geq \mu(B_2(\frac{1}{\sqrt{k}})),
$$
where $B_2(R)$ is the ball with radius $R$ for the norm $|.|_2$ in $H_1(\Gamma',\R)$. We deduce the inequality (2.4) for graphs.  \hspace{\stretch{1}} $\Box$\\

\bigskip

\noindent {\bf Acknowledgements}

\medskip

The author would like to thank his thesis advisor I.Babenko and express his gratitude to S.Sabourau for several helpful remarks.

{
}

\end{document}